\documentclass{article}

\usepackage{pinlabel}

\usepackage{graphicx}

\newcommand{\pic}[2]{\raisebox{-.5\height}{\includegraphics[scale=#2]{#1}}}

\def\Xor{\pic{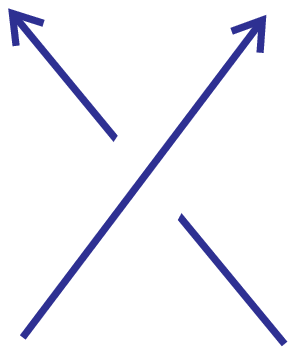}{.300}}
\def\Yor{\pic{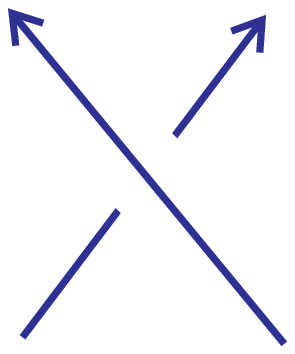} {.300}}
\def\Ior{\pic{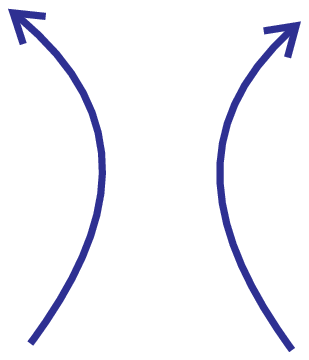} {.300}}

\def\Rcurlor{\pic{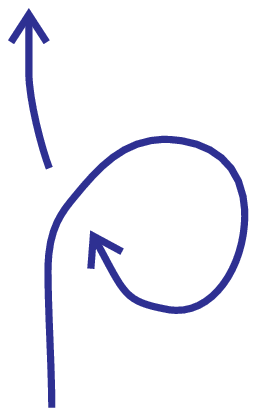}{.300}}
\def\Lcurlor{\pic{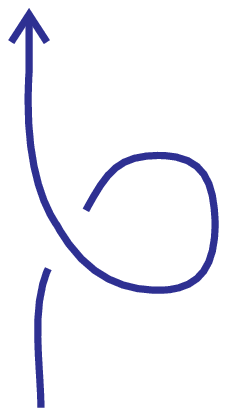} {.300}}
\def\Idor{\pic{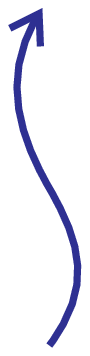} {.300}}

\def\Xunor{\pic{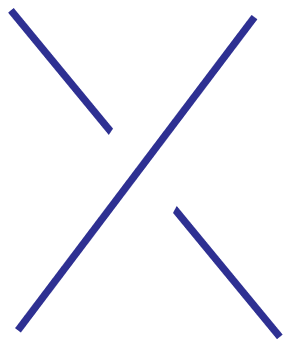}{.300}}
\def\Yunor{\pic{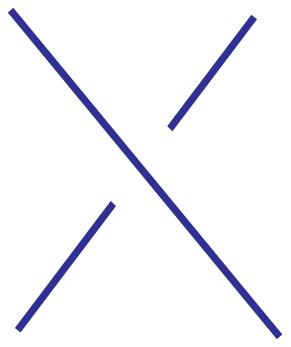} {.300}}
\def\Iunor{\pic{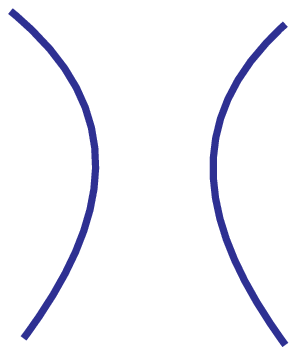} {.300}}
\def\Infunor{\pic{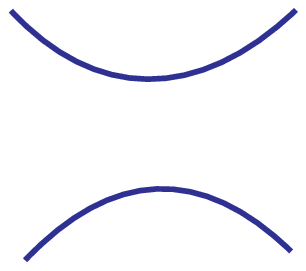} {.300}}

\def\Rcurlunor{\pic{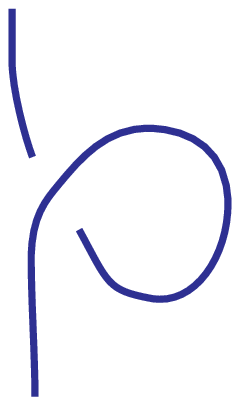}{.300}}
\def\Lcurlunor{\pic{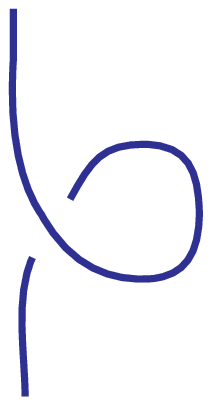} {.300}}
\def\Idunor{\pic{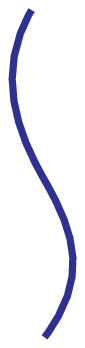} {.300}}

\def\Satellite{\pic{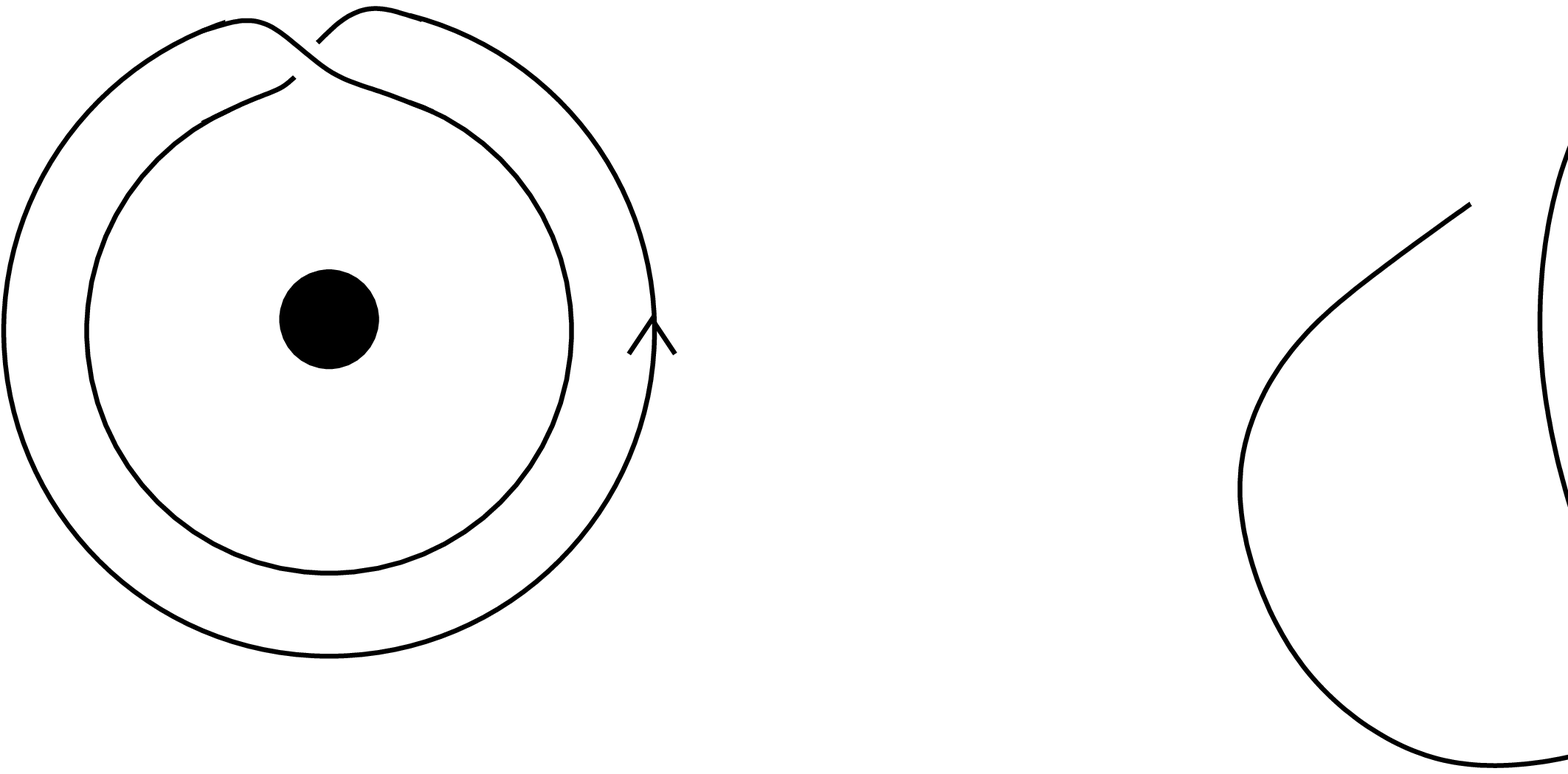}{.200}}

\def\Meridian{\pic{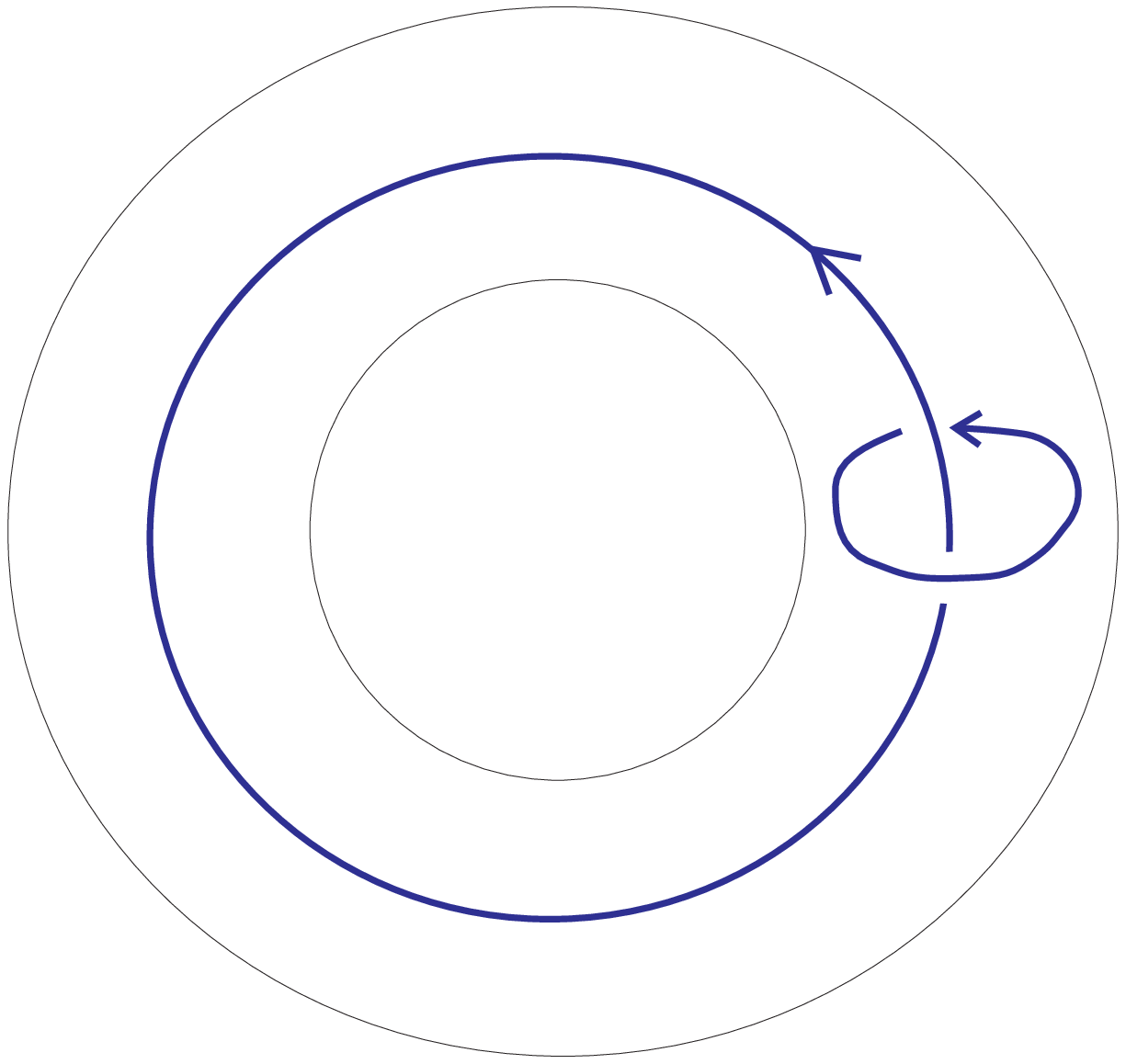}{.2500}}

\def\Meridianunor{\pic{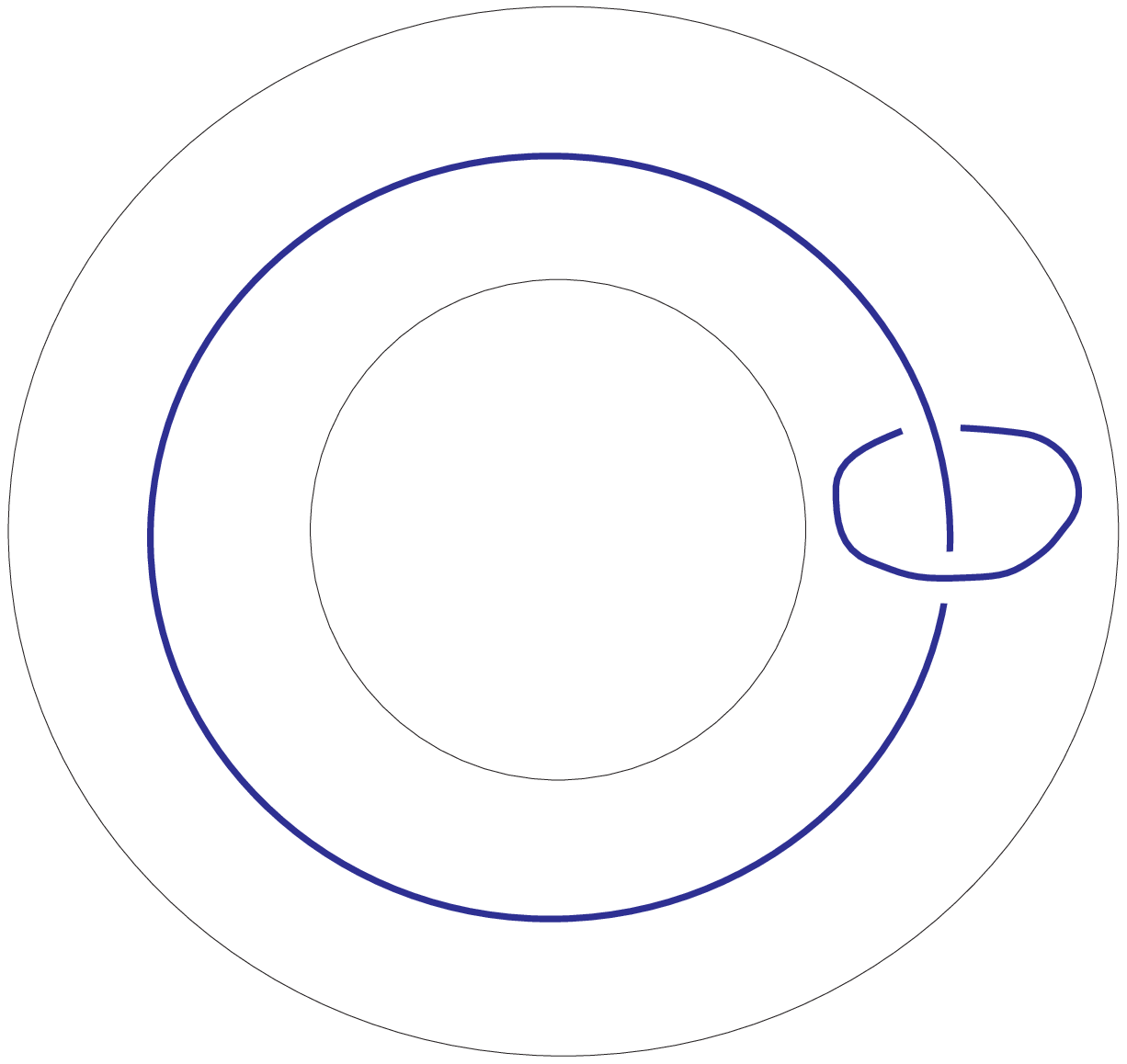}{.2500}}

\def\parallels{\pic{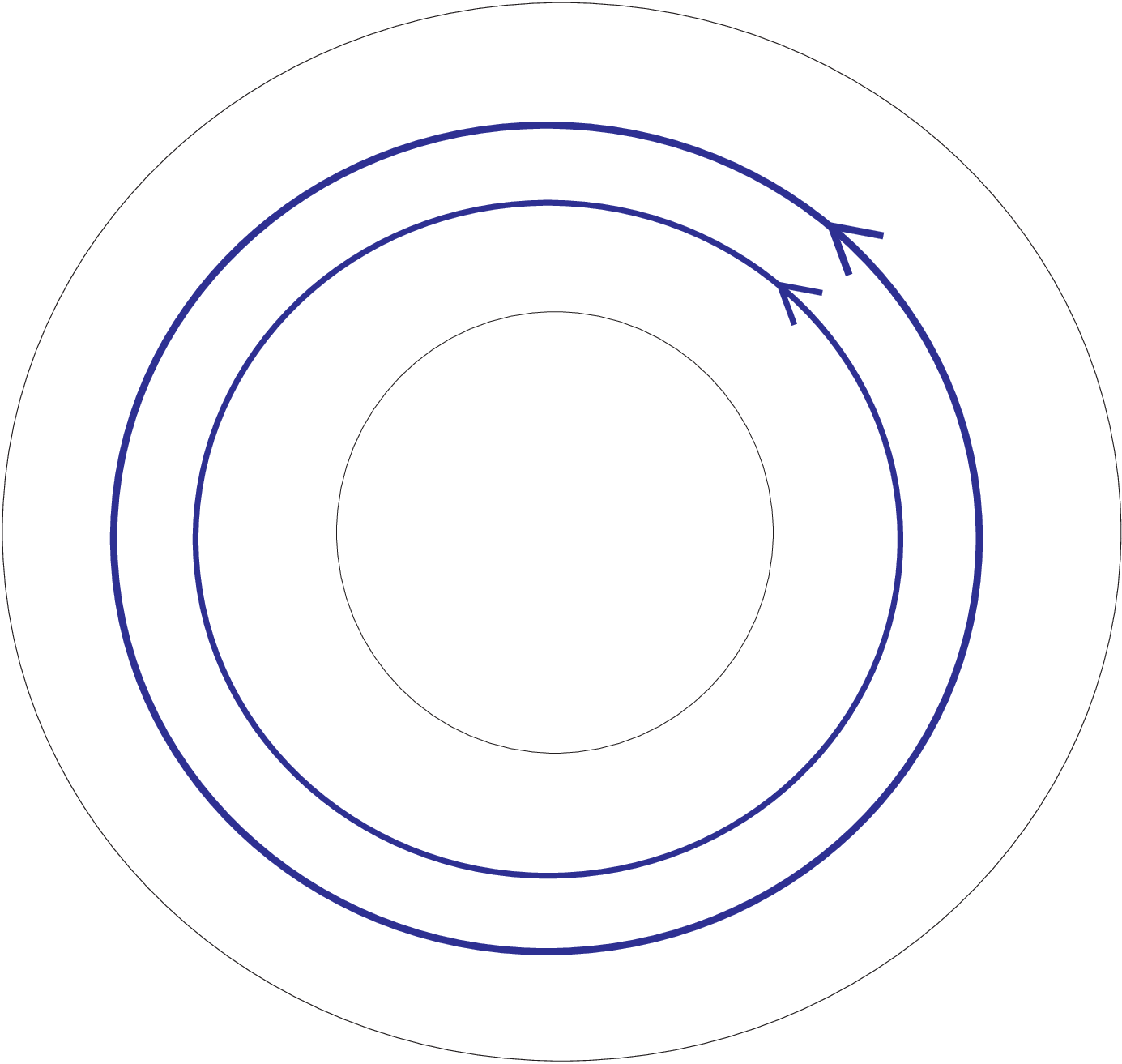}{.250}}

\def\LM{\pic{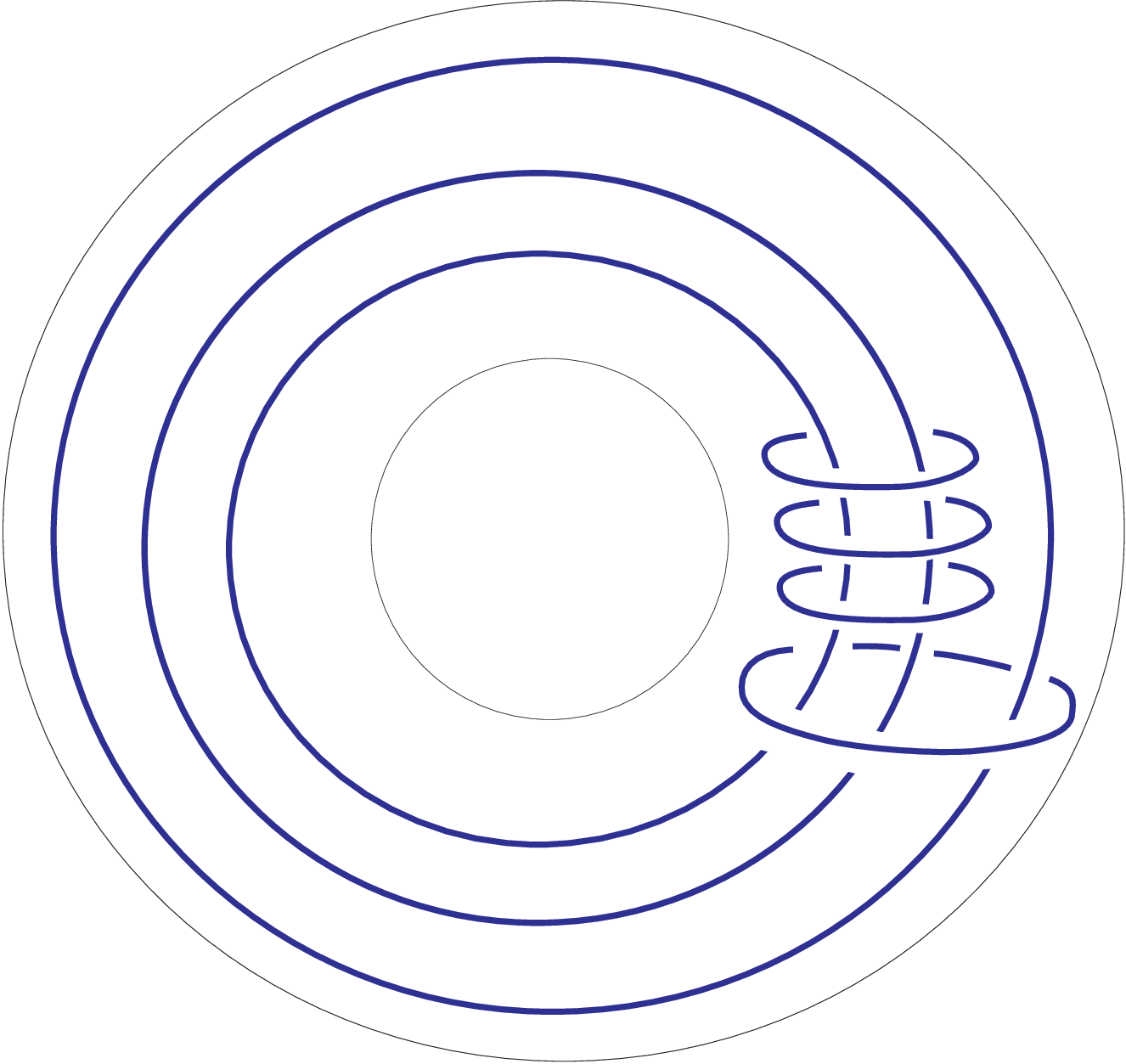}{.250}}
\def\lmthree{\pic{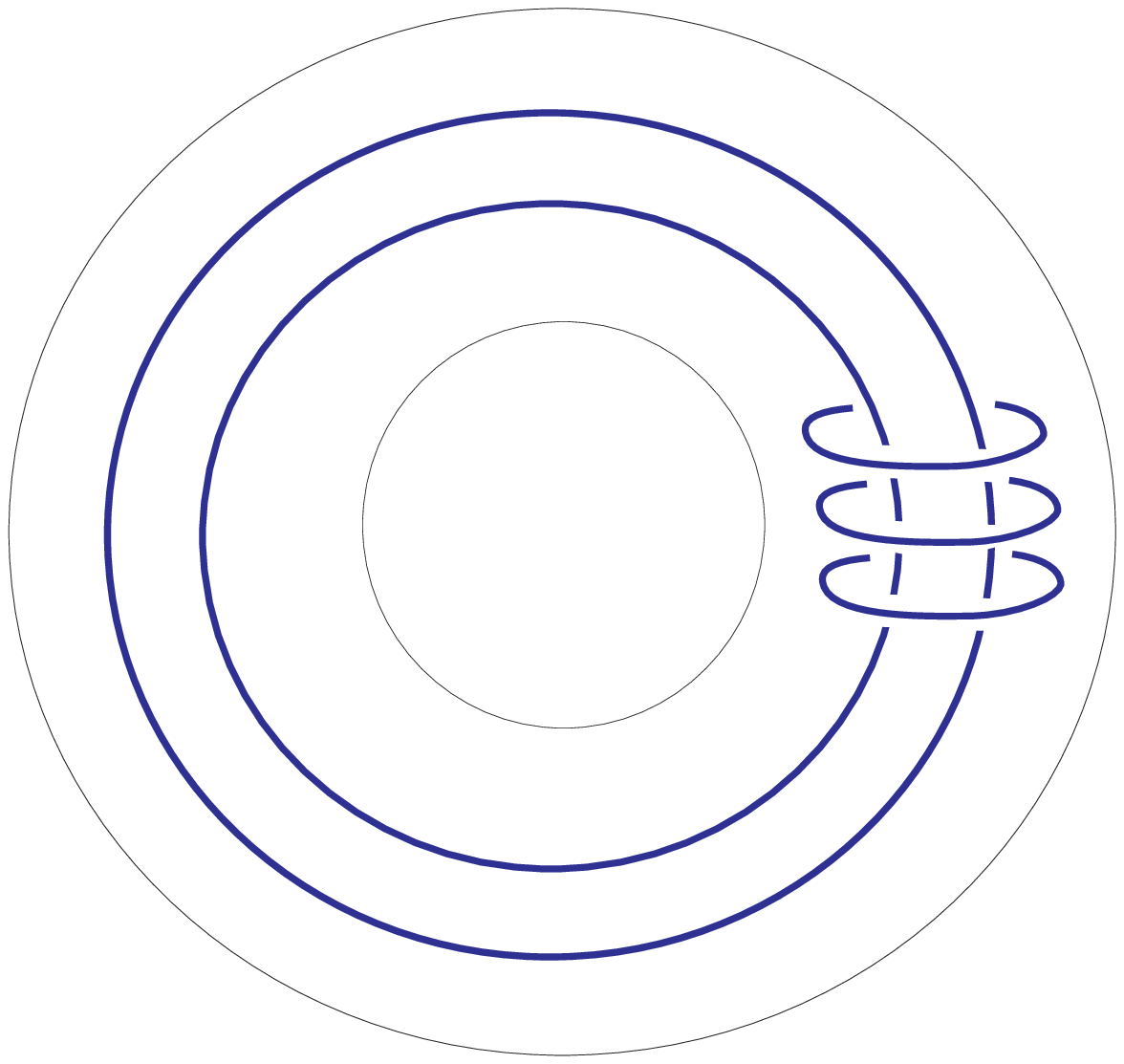}{.250}}

\newcommand{\bc}{\begin{center}}
\newcommand{\ec}{\end{center}}

\newcommand{\be}{\begin{equation}}
\newcommand{\ee}{\end{equation}}

\newcommand{\beqn}{\begin{eqnarray*}}
\newcommand{\eeqn}{\end{eqnarray*}}

\newcommand{\gl}{\lambda}

\newcommand{\ds}{\displaystyle}

\def\red{}
\def\blue{}
\def\redem{}

\usepackage{amsmath}
\usepackage{amsthm}

\newcommand{\qlm}{Q_{\gl,\mu}}

\newtheorem{theorem}{Theorem}

{\par\smallskip}
{\par\smallskip}
\renewenvironment{proof}[1][Proof]{\textit{#1.} }{\hfill \rule{0.5em}{0.5em}}

\begin{document}

\bc{\Large\bf Relations between Kauffman and Homfly satellite invariants\\[3mm]}
{\sc H. R. Morton {\rm and } N. D. A. Ryder\\[2mm]}
 {\small \sl Department of Mathematical Sciences\\ University of Liverpool\\ Peach Street, Liverpool L69 7ZL 
}
\ec
\begin{abstract}
We extend a mod $2$ relation between the Kauffman and Homfly polynomials, first observed by Rudolph in 1987, to the general Kauffman and Homfly \redem{satellite} invariants.
\end{abstract}

\section{Introduction}

The result of Rudolph \cite{Rudolph} relating the Kauffman polynomial of a link to the Homfly polynomials of reverse parallels of its sublinks has long provided an intriguing and rather isolated connection betwewen these two major knot invariants. In this paper we show how Rudolph's result can be extended to general Kauffman satellite invariants. 

Our results connect the meridian eigenvectors $y_\gl$ of the Kauffman skein of the annulus with the orientation-symmetric eigenvectors  $R_\gl$ of the Homfly skein in a combinatorially suggestive way, which begs for a representation-theoretic explanation. The proofs, however, rely purely on the representation-related combinatorics of these eigenvectors, and features of their eigenvalues, while depending eventually on Rudolph's original skein-theoretic result as a basis for the inductive proof.
 
  \subsection{Satellite invariants and skeins of the annulus}

We   use the framed version of the Homfly polynomial based on the skein relations
\bc $\Xor\quad -\quad\Yor \quad=\quad (s-s^{-1})\  \Ior$\\[2mm]
$\Rcurlor\quad = \quad v^{-1}\ \Idor\ , \qquad \Lcurlor\quad=\quad v\  \Idor\  $
\ec and normalised to take the value $1$ on the empty link.

Similarly the defining skein relations for the Kauffman polynomial are
\bc $\Xunor\quad -\quad\Yunor \quad=\quad (s-s^{-1}) \left(\Iunor-\Infunor\right)$\\[2mm]
$\Rcurlunor\quad = \quad v^{-1}\ \Idunor\ , \qquad \Lcurlunor\quad=\quad v\  \Idunor\  .$
\ec

Throughout the paper coefficients are taken in the ring   $\Lambda={\bf Z}_2[v^{\pm1},s^{\pm1}]$ with the elements $\{s^r-s^{-r},r\in {\bf Z}_{>0}\}$ inverted. As in Rudolph's original result the use of ${\bf Z}_2$ rather than ${\bf Z}$  is essential.

 A \emph{satellite} of a framed knot $K$ is determined by choosing a diagram $Q$ in the standard annulus, and then drawing $Q$ on the annular neighbourhood of $K$ determined by the framing. This gives the satellite of $K$ when \emph{decorated} by the diagram $Q$. The result is the satellite knot denoted by $K*Q$. 
{\labellist
\small
\pinlabel {$Q$} at -50 195
\pinlabel {$K$} at 440 195
\pinlabel {$K*Q$} at 1060 195
\endlabellist}
\bc
\Satellite
\ec

Starting with a framed link $$L=\cup_{i=1}^k L_i$$ we may decorate each component with a diagram $Q_i$ in the annulus. The Kauffman or Homfly polynomial of the resulting link, is called a \emph{satellite invariant} of $L$, as detailed for example in \cite{MortonHadji} in the case of the Homfly polynomial. We write $D(L;Q_1,\ldots,Q_k)$ and $P(L;Q_1,\ldots,Q_k)$ for the Kauffman and Homfly polynomials respectively of the link $L$ when each component $L_i$ is decorated by $Q_i$.

 The satellite invariants depend only on the decorations $Q_i$ as elements in the appropriate skein of the annulus (consisting of $\Lambda$-linear combinations of diagrams in the annulus, modulo the corresponding skein relations). The general linear combination of satellite invariants can be determined in terms of the satellite invariants where the decorations are drawn from  a linear basis of the skein of the annulus.

\subsection{The Kauffman skein}

The {Kauffman skein} of the annulus has a basis $\{{y_\gl}\}$, described by  Blanchet and Beliakova \cite{BlanchetBeliakova}, where $\gl$ runs through the set of all partitions. These elements consist of eigenvectors of the meridian map, \cite{Zhong}, and are closely related to the irreducible modules used in constructing link invariants from the $B,C$ and $D$ series of quantum groups.

When   $\gl$ is a partition of $1$ we write  $y_\gl=y_1$. This decoration can be represented in the skein of the annulus by  the single core curve of the annulus. Then $D(L;y_1,\ldots,y_1)$ is just the Kauffman polynomial of the link $L$ itself, with coefficients reduced mod 2.

\subsection{The Homfly skein}

The {Homfly skein} of the annulus has a basis $\{{Q_{\gl,\mu}}\}$, described by Hadji and Morton \cite{MortonHadji}, where $\gl$ and $\mu$ run through the set of all partitions, including the empty partition $\phi$. Diagrams in the Homfly skein are oriented. Reversing the orientation of all components in a diagram induces an involution on the skein which carries the basis element $\qlm$ to $Q_{\mu,\gl}$.

When   $\gl$ is a partition of $1$ and $\mu$ is empty the decoration $Q_{1,\phi}$ can be represented by a single oriented core curve, while $Q_{\phi,1}$ is represented by the oppositely oriented core curve, and $Q_{\phi,\phi}$ by the empty diagram.

The result presented here in Theorem \ref{Main} shows how the Kauffman invariant with decoration ${y_\gl}$ on a component of $L$ relates to the Homfly invariant with decoration $R_\gl=\red{Q_{\gl,\gl}}$ on the corresponding  component.

The decoration $R_1=Q_{1,1}$ for Homfly can be termed the \emph{adjoint decoration}, because of its relation with the adjoint representation for the $A$ series of quantum groups. In terms of the product in the skein of the annulus induced by placing decorations in concentric annuli we can write $Q_{1,1}=Q_{1,\phi}Q_{\phi,1}-Q_{\phi,\phi}$, which is realised as the decoration by two parallel oppositely oriented core curves minus the empty decoration.  The  \emph{adjoint Homfly polynomial}  $P_{\mbox{adj}}(L)$ is given by using the adjoint decoration $R_1$ on each component of an unoriented link $L$.

\subsection{Rudolph's relation and its extension}
Squaring in the ring $\Lambda$ is a ring homomorphism which we write as $\overline{\phantom {W}}:\Lambda\to\Lambda$.

Rudolph's original result \cite{Rudolph}  says that the adjoint Homfly polynomial of a link $L$   is equal to its Kauffman polynomial  in the ring $\Lambda$ after the substitution $v\to v^2, s\to s^2$. This can be summarised as $$P_{\mbox{adj}}(L)=\overline{D(L)}.$$

We extend Rudolph's result to general satellite invariants as follows.

\begin{theorem}\label{Main} For any choice of partitions $\gl(1),\ldots,\gl(k)$ the Homfly satellite invariant $P(L;R_{\gl(1)},\ldots,R_{\gl(k)})$ of the link $L$ satisfies

\beqn
P(L;R_{\gl(1)},\ldots,R_{\gl(k)})&=&\overline{D(L;y_{\gl(1)},\ldots,y_{\gl(k)})}\\
&=&D(L;y_{\gl(1)},\ldots,y_{\gl(k)})|_{v\to v^2,s\to s^2}.\eeqn
\end{theorem}

\subsection{The strategy for a proof of Theorem \ref{Main}}
The proof depends on  the very similar behaviour of the Kauffman and Homfly skeins of the annulus when a further meridian or longitude is adjoined. The effect is determined by the \redem{meridian maps}, established by Lu and Zhong \cite{Zhong} for Kauffman and by Hadji and Morton \cite{MortonHadji} for Homfly, and the \redem{branching rules}, described by Beliakova and Blanchet \cite{BlanchetBeliakova} for Kauffman and by Hadji \cite{Hadji} for Homfly.

 It is possible to produce a linear combination  $\sum a_j L^{(j)}$ of links $L^{(j)}$, each consisting of a succession of many longitude and meridian curves around each component $L_i$ of the original link $L$, with the property that 
$$D(L;y_{\gl(1)},\ldots,y_{\gl(k)})=\sum a_j D(L^{(j)})$$ and
$$P(L;R_{\gl(1)},\ldots,R_{\gl(k)})=\sum \overline{a_j} P_{\mbox{adj}}(L^{(j)}).$$ 
Rudolph's result applied to $L^{(j)}$ would then establish the theorem.

Here we prove Theorem \ref{Main} by induction on $\sum|\gl(i)-1|$, with Rudolph's result as the base, using an intermediate version of the longitude and meridian technique.

\section{Meridian maps}
The basic information needed about the skeins of the annulus are firstly the behaviour of the \redem{meridian maps}.

\subsection{Kauffman meridian maps} 
In the Kauffman skein Lu and Zhong \cite{Zhong} show that when $y_1$ is placed on a meridian around $y_\gl$ the result is $c_\gl y_\gl$, and the scalars $c_\gl$ for different partitions are all different.
\bc
{\labellist
\small
\pinlabel {$y_1$} at 372 595
\pinlabel {$y_\gl$} at 237 773
\endlabellist\Meridianunor}\quad = \quad $c_\gl y_\gl$.
\ec

  The eigenvalue $c_\gl$ calculated by Lu and Zhong uses  coefficients in ${\bf Z}[v^{\pm1},s^{\pm1}]$ with denominators as in $\Lambda$. 
\begin{theorem}[Lu and Zhong]\label{LZ}
\[
c_\gl = (s-s^{-1}) \left(v^{-1}\sum_{x \in \gl} s^{2c(x)} - v\sum_{x \in \gl} s^{-2c(x)} \right) + \frac{v^{-1}-v}{s-s^{-1}} + 1.\]
\end{theorem}
Here $x$ runs through cells in a diagram of the partition $\gl$, and $c(x)=j-i$ denotes the content of the cell $x$ in position $(i,j)$.

The result of Lu and Zhong  that these eigenvalues are distinct holds  even when they are evaluated with mod 2 coefficients in the ring $\Lambda$.
\begin{theorem} \label{distinct} If $c_\gl = c_\mu$ when evaluated in $\Lambda$ then $\gl=\mu$.
\end{theorem}
\begin{proof}  The \emph{content polynomial} $C_\gl(t)$ of a partition $\gl$ is the integer Laurent polynomial $$C_\gl(t)=\sum_{x\in\gl}t^{c(x)}.$$  As above, $c(x)$ denotes the content of the cell $x$ in a Young diagram of $\gl$. The eigenvalue $c_\gl$ can then be written as 
$$ c_\gl=(s-s^{-1})\left(v^{-1}C_\gl(s^2)-vC_\gl(s^{-2})\right)+ \frac{v^{-1}-v}{s-s^{-1}} + 1.$$ 

Write the partition $\gl$  in Frobenius form as   $(a_1,\ldots,a_k|b_1,\ldots,b_k)$, having arms of lengths $a_1>a_2>\ldots a_k\ge0$ and legs of lengths $b_1>b_2>\ldots b_k\ge0$.  Then $$(s-s^{-1})C_\gl(s^2)=\sum_{i=1}^k (s^{2a_i+1}-s^{-2b_i-1}),$$ since  the terms in the content polynomial $C_\gl(s^2)$ coming from the cells in the $i$th arm and leg of $\gl$ are $\ds\sum_{r=-b_i}^{a_i}s^{2r}$.

Then $(s-s^{-1})C_\gl(s^2)$ determines the Frobenius form of $\gl$, even when given as a Laurent polynomial in $s$  with ${\bf Z}_2$ coefficients, as there can be no cancellation among the $2k$ terms.

The coefficient of $v^{-1}$ in $c_\gl-c_\mu$ is $(s-s^{-1})C_\gl(s^2)-(s-s^{-1})C_\mu(s^2)$. If $c_\gl=c_\mu$ in $\Lambda$ then $(s-s^{-1})C_\gl(s^2)=(s-s^{-1})C_\mu(s^2)$ in ${\bf Z}_2[s^{\pm 1}]$, and so $\gl=\mu$.

\end{proof}
\subsection{Homfly meridian maps}
In the Homfly skein Hadji and Morton \cite{MortonHadji} show that the basis element $\qlm$ is an eigenvector for the meridian map on the Homfly skein of the annulus induced by placing an oriented meridian around the annulus.

\bc
{\labellist
\small 
\pinlabel {$\qlm$} at 237 773
\endlabellist\Meridian}\quad = \quad $s_{\gl,\mu}\qlm$.

\ec
\begin{theorem}[Hadji and Morton]\label{HM}
The eigenvalue $s_{\gl,\mu}$ is given  by the formula
\[
s_{\gl,\mu} = (s-s^{-1}) \left(v^{-1}C_\gl(s^2)  - v C_\mu(s^{-2})  \right) + \frac{v^{-1}-v}{s-s^{-1}}.\]

Orienting the meridian in the opposite sense results in the eigenvalue $s_{\mu,\gl}$.
\end{theorem}
Comparing the formulae in Theorems \ref{LZ} and \ref{HM} shows immediately that $c_\gl=s_{\gl,\gl}+1$.

\begin{theorem} When the adjoint decoration $R_1$ is placed on a meridian around $\qlm$ the resulting element in the Homfly skein is $(s_{\gl,\mu}s_{\mu,\gl}-1)\qlm$.
\end{theorem}
\begin{proof} Since $R_1=Q_{1,\phi}Q_{\phi,1}-1$ we can use the eigenvalue results above.
\end{proof}
 
\begin{theorem}
 When the adjoint decoration $R_1$ is placed on a meridian around $R_\gl$ the resulting element  in the Homfly skein  using the ring $\Lambda$ is $\overline{c_\gl}R_\gl$.
\bc
{\labellist
\small
\pinlabel {$R_1$} at 372 595
\pinlabel {$R_\gl$} at 237 773
\endlabellist\Meridian}\quad = \quad $\overline{c_\gl} R_\gl$.

\ec
\end{theorem}

\begin{proof} Since $R_\gl=Q_{\gl,\gl}$ we have 
\bc
{\labellist
\small
\pinlabel {$R_1$} at 372 595
\pinlabel {$R_\gl$} at 237 773
\endlabellist\Meridianunor}\quad = \quad $(s_{\gl,\gl}^2 -1) R_\gl$.
\ec
In the ring $\Lambda$ we have $\overline{c_\gl}=c_\gl^2=s_{\gl,\gl}^2+1=s_{\gl,\gl}^2 -1$.
\end{proof}

\section{Branching rules}
  
The skeins of the annulus admits a commutative product, induced by placing two skein elements on parallel strands.
When a single {parallel} strand is added in the skein the description of the resulting product is determined by the \redem{branching rules}.

\subsection{Kauffman branching rules}

\begin{figure}[ht]
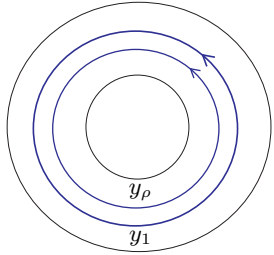

\bc
{\labellist
\small
\pinlabel {$y_1$} at 223 458
\pinlabel {$y_\rho$} at 223 528
\endlabellist\parallels}\ec
\caption{The product $y_\rho y_1$}\label{kauffmanbranch}
\end{figure}

In the {Kauffman} skein of the annulus two parallel strands, one decorated by $\red{y_\rho}$ and the other by ${y_1}$ as in Figure \ref{kauffmanbranch}, can be written as a sum of basis elements using the branching rules of Beliakova and Blanchet \cite{BlanchetBeliakova}.

 The branching formula for the product $y_\rho y_1$ is given in the following theorem.

\begin{theorem}[Beliakova and Blanchet]
$$y_\rho y_1=\sum_{\mu\in\rho^+\cup\rho^-}y_\mu,$$ where $\rho^+$ and $\rho^-$ are the sets of partitions given from $\rho$ by respectively adding  or deleting one cell.
\end{theorem}

\subsection{Homfly branching rules}

A single oriented core curve in the Homfly skein of the annulus represents $Q_{1,\phi}$ or $Q_{\phi,1}$ depending on its direction.

\begin{figure}[ht]
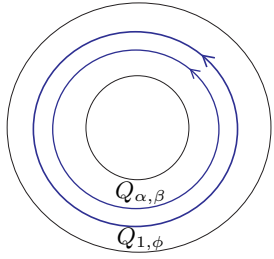

\bc
{\labellist
\small
\pinlabel {$Q_{1,\phi}$} at 223 458
\pinlabel {$Q_{\alpha,\beta}$} at 223 528
\endlabellist\parallels}\ec
\caption{The product $Q_{\alpha,\beta}Q_{1,\phi}$}\label{Homflybranch}
\end{figure}

The basic branching rules  
in the {Homfly} skein for two parallel strands, one decorated by $Q_{\alpha,\beta}$ and one by $Q_{1,\phi}$, as in Figure \ref{Homflybranch}, are  
$$Q_{\alpha,\beta}Q_{1,\phi}=\sum_{\mu\in\alpha^+ }Q_{\mu,\beta}+\sum_{\nu\in\beta^-}Q_{\alpha,\nu}, $$ while reversing the direction in the annulus gives a similar formula for the product with $Q_{\phi,1}$. 

We can then expand the product $R_\rho R_1=Q_{\rho,\rho}(Q_{1,\phi}Q_{\phi,1}-1)$ as follows.
  
\begin{theorem}
\label{thm:hsr}
\[R_\rho R_1 =   \sum_{\mu \in \rho^+ \cup \rho^- } R_\mu+ \sum_{\substack{\{\alpha,\beta\} \\ \alpha \ne \beta}} n_{\alpha,\beta}(Q_{\alpha,\beta} + Q_{\beta,\alpha}) + 2|\rho^-|R_\rho, \] for some integers $n_{\alpha,\beta}$.
\end{theorem}

The main importance of this expansion for us is its similarity to the expansion for $y_\rho y_1$ in the Kauffman skein and the fact that the remaining terms occur in pairs $Q_{\alpha,\beta}+Q_{\beta,\alpha}$. There is an explicit calculation in \cite{Ryder} of the coefficients $n_{\alpha,\beta}$, which are either $1$ or $0$, but the details are not needed here. 

\section{Longitude-meridian decorations}

Using the branching rules and the meridian map it is possible to express any element in the Kauffman skein of the annulus as a $\Lambda$-linear combination of \emph{longitude-meridian} diagrams.

\begin{figure}[ht]
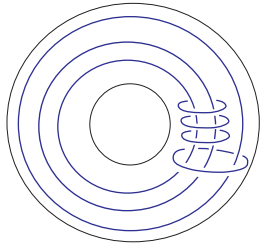

\bc
\LM
\ec
\caption{The longitude-meridian diagram $l^2m^3lm$}\label{LMfig}
\end{figure}

We define a longitude-meridian diagram recursively, starting from a single core curve, by including a sequence of further longitude and meridian curves. These can be specified by a sequence of powers of $l$ and $m$, as in Figure \ref{LMfig}.

For a partition $\gl$ with $|\gl|>1$ choose a subpartition $\rho$ with one fewer cell. We can reproduce the effect of the single element $y_\gl$ in the Kauffman skein of the annulus in terms of a linear combination of longitude-meridian diagrams in which  each component is decorated by $y_1$ except for one innermost longitude decorated by $y_\rho$ as follows.

  Write \be \label{eq1} X(t)=\prod_{\mu\in\rho^+\cup\rho^--\{\gl\}}(t-c_\mu)\ee as a polynomial $\sum \blue{a_r} t^r$ in $t$. 

The longitude-meridian diagram $l^2m^r$  decorated with $y_\rho$ on one longitude and with $y_1$ on the other curves is the product $y_\rho y_1$ surrounded by $r$ meridians.
 Since $\ds y_\rho y_1=\sum_{\mu\in\rho^+\cup\rho^-} y_\mu$ this decorated longitude-meridian diagram represents $\ds \sum_{\mu\in\rho^+\cup\rho^-} c_\mu^r y_\mu$ in the Kauffman skein of the annulus.

Then the linear combination $\ds\sum_r a_r l^2 m^r$ of longitude-meridian diagrams decorated with $y_\rho $ on one longitude and with $y_1$ elsewhere represents
\beqn \sum_r a_r \left(\sum_{\mu\in\rho^+\cup\rho^-} c_\mu^r y_\mu\right)&=&\sum_{\mu\in\rho^+\cup\rho^-} X(c_\mu) y_\mu \\&=&X(c_\gl)y_\gl, \mbox{ by }(\ref{eq1}).\eeqn

Since $X(c_\gl)\ne 0$ in $\Lambda$, by Theorem \ref{distinct}, this represents $y_\gl$ up to a non-zero multiple. 

\begin{figure}[ht]
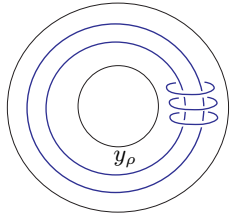
\label{lm3}
\bc{\labellist
\small
\pinlabel {$y_\rho$} at 223 540
\endlabellist
\lmthree}
\ec
\caption{The decoration for the coefficient of $t^3$ in $X(t)$.}
\end{figure}

For example, the decoration for $r=3$ is shown in Figure \ref{lm3}.

Although we will not need to expand further it is clear that    the longitude decorated by $y_\rho$ could itself be replaced recursively by a linear combination of longitude-meridian diagrams, to give a representation of $y_\gl$ in the skein entirely in terms of longitude-meridian diagrams.

\section{The proof of Theorem \ref{Main}}
We proceed by induction on $N=\sum(|\gl(i)|-1)$.  Components which are assigned the partition $1$ do not then contribute to $N$.  The base case with $N=0$, where all components have the partition $1$, is given by Rudolph's result.

If $N>0$ then $|\gl(i)|>1$ for some $i$. We will suppose that $i=1$ and  write $\gl=\gl(1)$. 

  Write $L^{(r)}$ for the link consisting of two parallel copies of $L_1$ surrounded by $r$ meridians, along with the remaining components of $L$. This is the satellite of $L$ in which the component $L_1$ is decorated by the longitude-meridian diagram $l^2m^r$.

Choose a partition $\rho\in \gl^-$. This is possible since $|\gl|>1$.
Decorate $L^{(r)}$   with $y_\rho$ on one parallel of $L_1$, $y_1$ on the other parallel and on the $r$ meridians, and with the original decorations $y_{\gl(j)}$ on $L_2,\ldots, L_k$, and write $D^{(r)}$ for its Kauffman invariant.

Write $P^{(r)}$ for the Homfly polynomial of $L^{(r)}$ when decorated by the corresponding elements $R_\rho$, $R_1$ and $R_{\gl(j)}$.

The element in the Homfly skein represented by the diagram $l^2m^r$ decorated by $R_\rho$ on one longitude and $R_1$ on all the other curves represents
\[ \sum_{\mu \in \rho^+ \cup \rho^- }(s_{\mu,\mu}^2-1)^r R_\mu+ \sum_{\substack{\{\alpha,\beta\} \\ \alpha \ne \beta}} n_{\alpha,\beta}(s_{\alpha,\beta}s_{\beta,\alpha}-1)^r(Q_{\alpha,\beta} + Q_{\beta,\alpha}) + 2|\rho^-|(s_{\rho,\rho}^2-1)^r R_\rho .\] 
Hence we can write $P^{(r)}$  as the sum of invariants
\beqn
P^{(r)}&=&\sum_{\mu \in \rho^+ \cup \rho^- }(s_{\mu,\mu}^2-1)^r P(L;R_\mu,R_{\gl(2)},\ldots,R_{\gl(k)})\\
&&+ \sum_{\substack{\{\alpha,\beta\} \\ \alpha \ne \beta}} n_{\alpha,\beta}(s_{\alpha,\beta}s_{\beta,\alpha}-1)^r(P(L;Q_{\alpha,\beta},R_{\gl(2)},\ldots,R_{\gl(k)}) \\
&&+ P(L;Q_{\beta,\alpha},R_{\gl(2)},\ldots,R_{\gl(k)})) \\&&+ 2|\rho^-|(s_{\rho,\rho}^2-1)^r P(L;R_\rho,R_{\gl(2)},\ldots,R_{\gl(k)}).
\eeqn

Since the Homfly polynomial of a link is unchanged when the orientation of \emph{all} its components are reversed we have $$P(L;Q_{\alpha,\beta},R_{\gl(2)},\ldots,R_{\gl(k)})=P(L;Q_{\beta,\alpha},R_{\gl(2)},\ldots,R_{\gl(k)}).$$ Then in $\Lambda$ we have the simplified expression \be\label{eq2} P^{(r)}=\sum_{\mu \in \rho^+ \cup \rho^- }(\overline{c_\mu})^r P(L;R_\mu,R_{\gl(2)},\ldots,R_{\gl(k)}),\ee
since $s_{\mu,\mu}^2-1=\overline{c_\mu}$ and the remaining terms cancel mod 2.

  The induction count $N$  has reduced by $1$ for the choice of partitions for the links $L^{(r)}$ when compared with the choice for the original link $L$.

Then $P^{(r)}=\overline{D^{(r)}}$, by induction.

Define the polynomial  $$X(t)=\sum a_r t^r=\prod_{\mu\in\rho^+\cup\rho^--\{\gl\}}(t-c_\mu),$$ as in equation (\ref{eq1}). Then
  $$\sum {a_r} D^{(r)}=X(c_\gl)D(L;y_{\gl },y_{\gl(2)},\ldots,y_{\gl(k)}).$$
Now \beqn \sum \overline{ a_r D^{(r)}}&=&\sum \overline{ a_r} \overline{D^{(r)}}\\
&=&\sum \overline{a_r} P^{(r)}\\
&=&\sum_{\mu \in \rho^+ \cup \rho^- }(\sum\overline{a_r}(\overline{c_\mu})^r )P(L;R_\mu,R_{\gl(2)},\ldots,R_{\gl(k)}), \mbox{ by } (\ref{eq2}),\\
&=&\sum_{\mu \in \rho^+ \cup \rho^- }\overline{X(c_\mu)}P(L;R_\mu,R_{\gl(2)},\ldots,R_{\gl(k)}) \\
&=&\overline{X(c_\gl)}P(L;R_\gl,R_{\gl(2)},\ldots,R_{\gl(k)}), \mbox{ by }(\ref{eq1}).
\eeqn
Then $\overline{X(c_\gl)}\ \overline{D(L;y_{\gl },y_{\gl(2)},\ldots,y_{\gl(k)})}=\overline{X(c_\gl)}P(L;R_\gl,R_{\gl(2)},\ldots,R_{\gl(k)})$.

Theorem \ref{Main} follows immediately, since $X(c_\gl)\ne 0$ in $\Lambda$, and $\Lambda$ has no zero-divisors.
\section*{Acknowledgments}

The thought of using a longitude and meridian construction to realise the effect of a prescribed decoration in either skein of the annulus was triggered by a talk of Vogel. Coupled with the close similarity observed between the eigenvalues for $y_\gl$ and $R_\gl$ under the respective meridian maps, as calculated in \cite{Zhong} and \cite{MortonHadji}, this led the first author  to the conjectured extension \cite{integrality} of Rudolph's result, proved here in Theorem \ref{Main}.

The detailed work and final confirmation of the conjecture form a part of the second author's PhD thesis. Work for this was supported by EPSRC funding through Liverpool University under the doctoral training grant number EP/P500338/1.

\end{document}